\documentclass[12pt]{amsart}
\usepackage{epic,eepic}
\def\draft{n}

\newtheorem{prop}{Proposition}[section]
\newtheorem{theorem}[prop]{Theorem}

\newtheorem{cor}[prop]{Corollary}
\newtheorem{definition}[prop]{Definition}

\def\lbl#1{\label{#1}\printname{#1}}

\def\printname#1{
        \if\draft y
                \smash{\makebox[0pt]{\hspace{-0.5in}
                        \raisebox{8pt}{\tt\tiny #1}}}
        \fi
}

\newcommand{\mathmode}[1]{$#1$}
\newlength{\standardunitlength}
\catcode`\@=11
\font\bfcaptionfont=cmssbx10 scaled \magstephalf
\newdimen\@captionmargin\@captionmargin=2\parindent
\long\def\@makecaption#1#2{
    \vskip 10pt
    \setbox\@tempboxa\hbox{%
      \small\sf{\bfcaptionfont #1. }\ignorespaces #2}
    \ifdim \wd\@tempboxa >\captionwidth {
        \rightskip=\@captionmargin\leftskip=\@captionmargin
        \unhbox\@tempboxa\par}
      \else
        \hbox to\hsize{\hfil\box\@tempboxa\hfil}
    \fi}
\catcode`\@=12

\advance \topmargin by -\headheight
\advance \topmargin by -\headsep
\evensidemargin \oddsidemargin
\begin{document}
\newdimen\captionwidth\captionwidth=\hsize

\title
  [Finite Type Invariants of Links with a Fixed Linking Matrix]
  {Finite Type Invariants\\of Links with a Fixed Linking Matrix}
\author{Eli Appleboim}
\address{Department of Mathematics\\
  Technion, Haifa 31000\\
  Israel
}
\email{eaea@techunix.technion.ac.il}

\thanks{This paper is available electronically at
  {\tt http://xxx.lanl.gov/abs/math/9906138}.
}
\thanks{Submitted to the Journal of Knot Theory and its Ramifications on June99}

\date{This edition: \today; }

\begin{abstract}
In this paper we introduce two theories of finite type invariants for framed 
links
with a fixed linking matrix. We show that these theories are different
from, but related to, the theory of Vassiliev invariants of knots and
links. We will take special note of the case of zero linking matrix.
i.e., zero-framed algebraically split links. We also study the
corresponding spaces of ``chord diagrams''.
\end{abstract}

\maketitle

\tableofcontents

\section{Introduction}

We divide the set of all framed link types (isotopy classes of framed
links in $S^{3}$) into equivalence classes according to their linking
matrices. Namely, every equivalence class consists of all framed link
types having a certain given linking matrix.

There exist for sometime now, definitions of finite type invariants of
topological objects such as knots, 3-manifolds, graphs etc. In order to
make these definitions, one should have a notion of adjacency. Then,
based on this notion it is possible to set conditions for an
invariant of a certain type of objects (knots, 3-manifolds etc.) to
be of finite type. Recall the standard definition for links (see
e.g.~\cite{ba:1, Birman:Bulletin, BirmanLin:Vassiliev, Goussarov:New,
Goussarov:nEquivalence, Kontsevich:Vassiliev, Vassiliev:CohKnot,
Vassiliev:Book} and~\cite{Bar-Natan:VasBib}). Two links are adjacent
if they differ by a single crossing change. The definition of Vassiliev
invariants of links is based on this notion of adjacency. However, if one
wishes to set up a finite type theory for links within a given linking
class, this notion of adjacency will not do since one crossing change
changes some linking number and therefore alters the linking matrix.

In this paper we introduce two forms of finite type invariants of links
in a fixed linking class, based on similar, yet slightly different,
notions of adjacency.

As a special case we will be interested in zero-framed algebraically
split links (zero framed links having zero linking matrix). Such links
appear in a very essential form in the theory of finite type invariants
of $3$-manifolds.\footnote{To be precise, the links that appear in
the theory of finite type invariants of 3-manifolds are unit framed.
The difference between unit-framed and zero-framed plays no role in the
statements we make in this paper.}

\begin{figure}
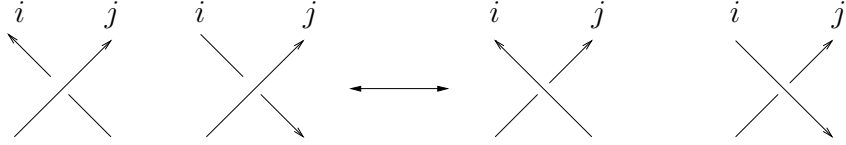

\[ 
  \printname{figure11}
  \setlength{\unitlength}{0.5\standardunitlength}
  \begin{array}{c}  \hspace{-1.7mm}
    \raisebox{-8pt}{\input draws/figure11.tex }
    \hspace{-1.9mm}
  \end{array}
 \]
\caption{A DD-crossing change} \label{dd-crossing}
\end{figure}

\begin{figure}
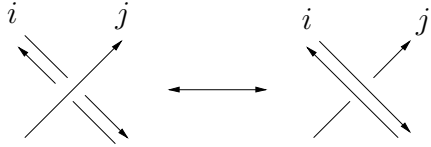

\[ 
  \printname{figure12}
  \setlength{\unitlength}{0.5\standardunitlength}
  \begin{array}{c}  \hspace{-1.7mm}
    \raisebox{-8pt}{\input draws/figure12.tex }
    \hspace{-1.9mm}
  \end{array}
 \]
\caption{A $\bigwedge$-crossing change} \label{w-crossing}
\end{figure}

\begin{definition}
\begin{enumerate}
\item  We will say that two links in the same linking class are
DD-adjacent if they differ by a double crossing change, as indicated
in Figure~\ref{dd-crossing}, showing two crossings involving the
$i^{th}$ and $j^{th}$ components in such a way that their linking
number remains unchanged.

\item The other type of the adjacency we present here is
influenced by the way the crossings to be changed are placed
along the link. We will call the two crossings involved ``close'' if
there is a segment of the $i^{th}$ component, say, going from one
of the crossings to the other one, meeting no other crossings along the
way.  This type of crossing change can be visualized as in
Figure~\ref{w-crossing}. In this case the links are called $\bigwedge$-adjacent.
\end{enumerate}
\end{definition}

Notice that since on the $i^{th}$ component the crossings are
close, we can think of them as appearing almost at the same
point. This will be of significance while considering the
corresponding chord diagrams in this case.

As an intermediate step, for both kinds of crossing changes, we get a
link with two double points. Such a link is the basic element of the
singular objects in the theory. Still, in order to be able to extend
link invariants to singular objects and to set a finiteness condition
we need to be more precise in the definition of a singular link in a
linking class.  Let ${\mathcal P}$ denote the class of links with some
fixed linking matrix $P$.

\begin{definition}
\begin{enumerate}
\item \label{dd-sing}
A $DD^{m}$-singular P-link is a link having
$2m$ double points paired
into $m$ ordered (by $+/-$ signs) pairs, so that each pair is of the form
\par \[ 
  \printname{figure1}
  \setlength{\unitlength}{0.5\standardunitlength}
  \begin{array}{c}  \hspace{-1.7mm}
    \raisebox{-8pt}{%
\begingroup\makeatletter\ifx\SetFigFont\undefined
\def\x#1#2#3#4#5#6#7\relax{\def\x{#1#2#3#4#5#6}}%
\expandafter\x\fmtname xxxxxx\relax \def\y{splain}%
\ifx\x\y   
\gdef\SetFigFont#1#2#3{%
  \ifnum #1<17\tiny\else \ifnum #1<20\small\else
  \ifnum #1<24\normalsize\else \ifnum #1<29\large\else
  \ifnum #1<34\Large\else \ifnum #1<41\LARGE\else
     \huge\fi\fi\fi\fi\fi\fi
  \csname #3\endcsname}%
\else
\gdef\SetFigFont#1#2#3{\begingroup
  \count@#1\relax \ifnum 25<\count@\count@25\fi
  \def\x{\endgroup\@setsize\SetFigFont{#2pt}}%
  \expandafter\x
    \csname \romannumeral\the\count@ pt\expandafter\endcsname
    \csname @\romannumeral\the\count@ pt\endcsname
  \csname #3\endcsname}%
\fi
\fi\endgroup
{\renewcommand{\dashlinestretch}{30}
\begin{picture}(2493,944)(0,-10)
\put(609,383){\blacken\ellipse{70}{70}}
\put(609,383){\ellipse{70}{70}}
\put(1809,383){\blacken\ellipse{70}{70}}
\put(1809,383){\ellipse{70}{70}}
\path(1509,683)(2109,83)
\path(2002.934,146.640)(2109.000,83.000)(2045.360,189.066)
\path(309,83)(909,683)
\path(845.360,576.934)(909.000,683.000)(802.934,619.360)
\path(909,83)(309,683)
\path(415.066,619.360)(309.000,683.000)(372.640,576.934)
\path(1509,83)(2109,683)
\path(2045.360,576.934)(2109.000,683.000)(2002.934,619.360)
\put(1209,83){\makebox(0,0)[lb]{\smash{{\mathmode{,}}}}}
\put(834,308){\makebox(0,0)[lb]{\smash{{\mathmode{+}}}}}
\put(1959,308){\makebox(0,0)[lb]{\smash{{\mathmode{-}}}}}
\put(309,833){\makebox(0,0)[lb]{\smash{{\mathmode{i}}}}}
\put(909,833){\makebox(0,0)[lb]{\smash{{\mathmode{j}}}}}
\put(1434,833){\makebox(0,0)[lb]{\smash{{\mathmode{i}}}}}
\put(2109,833){\makebox(0,0)[lb]{\smash{{\mathmode{j}}}}}
\put(1734.000,458.000){\arc{1500.000}{5.6397}{6.9267}}
\put(759.000,458.000){\arc{1500.000}{2.4981}{3.7851}}
\end{picture}
}
 }
    \hspace{-1.9mm}
  \end{array}
 \]
where $i$ and $j$ are components of the link, and so that
resolving those pairs in the way described below ends with a link
in the class ${\mathcal P}$.
\item A $\bigwedge^{m}$-singular P-link is the same as in~\ref{dd-sing}
but with the exception that all the pairs of double points taken in account
are close pairs, as shown in the figure below.
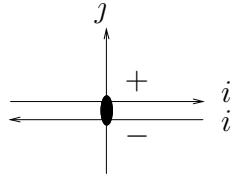
\[ 
  \printname{figure1a}
  \setlength{\unitlength}{0.5\standardunitlength}
  \begin{array}{c}  \hspace{-1.7mm}
    \raisebox{-8pt}{%
\begingroup\makeatletter\ifx\SetFigFont\undefined
\def\x#1#2#3#4#5#6#7\relax{\def\x{#1#2#3#4#5#6}}%
\expandafter\x\fmtname xxxxxx\relax \def\y{splain}%
\ifx\x\y   
\gdef\SetFigFont#1#2#3{%
  \ifnum #1<17\tiny\else \ifnum #1<20\small\else
  \ifnum #1<24\normalsize\else \ifnum #1<29\large\else
  \ifnum #1<34\Large\else \ifnum #1<41\LARGE\else
     \huge\fi\fi\fi\fi\fi\fi
  \csname #3\endcsname}%
\else
\gdef\SetFigFont#1#2#3{\begingroup
  \count@#1\relax \ifnum 25<\count@\count@25\fi
  \def\x{\endgroup\@setsize\SetFigFont{#2pt}}%
  \expandafter\x
    \csname \romannumeral\the\count@ pt\expandafter\endcsname
    \csname @\romannumeral\the\count@ pt\endcsname
  \csname #3\endcsname}%
\fi
\fi\endgroup
{\renewcommand{\dashlinestretch}{30}
\begin{picture}(2683,2085)(0,-10)
\put(1212,800){\blacken\ellipse{150}{374}}
\put(1212,800){\ellipse{150}{374}}
\put(2637,912){\makebox(0,0)[lb]{\smash{{\mathmode{i}}}}}
\put(2637,537){\makebox(0,0)[lb]{\smash{{\mathmode{i}}}}}
\put(1062,1962){\makebox(0,0)[lb]{\smash{{\mathmode{j}}}}}
\path(2412,687)(12,687)
\path(132.000,717.000)(12.000,687.000)(132.000,657.000)
\path(1212,12)(1212,1812)
\path(1242.000,1692.000)(1212.000,1812.000)(1182.000,1692.000)
\path(12,912)(2412,912)
\path(2292.000,882.000)(2412.000,912.000)(2292.000,942.000)
\put(1437,1062){\makebox(0,0)[lb]{\smash{{\mathmode{+}}}}}
\put(1437,387){\makebox(0,0)[lb]{\smash{{\mathmode{-}}}}}
\end{picture}
}
 }
    \hspace{-1.9mm}
  \end{array}
 \]
\end{enumerate}
Resolving a pair of singular points is done, in both types of singularity by,
\begin{enumerate}
\item
For DD-singularity,
\par \[ 
  \printname{figure2}
  \setlength{\unitlength}{0.5\standardunitlength}
  \begin{array}{c}  \hspace{-1.7mm}
    \raisebox{-8pt}{\input draws/figure2.tex }
    \hspace{-1.9mm}
  \end{array}
 \]
\item
For $\bigwedge$-singularity,
\par \[ 
  \printname{figure03}
  \setlength{\unitlength}{0.5\standardunitlength}
  \begin{array}{c}  \hspace{-1.7mm}
    \raisebox{-8pt}{\input draws/figure03.tex }
    \hspace{-1.9mm}
  \end{array}
 \]
\end{enumerate}
\end{definition}
\begin{definition} \label{def-3}
If $V$ is a $G$-valued link invariant where $G$ is any commutative
group then we extend $V$ to singular links, in both types by,
\begin{equation}
V(L) = V(\overline L) - V(\underline L)
\end{equation}
\end{definition}
\begin{definition} \label{fi-type}
We say that $V$ is a DD- linking class-P invariant of type
$m$ ($\bigwedge$-invariant of type $m$) if it vanishes on all
$DD^{m}$-singular P-links (all $\bigwedge^{m}$-singular P-links) We say
that $V$ is a DD($\bigwedge$)-linking class-P invariant of finite type
if it is of type $m$ for some positive $m$.
\end{definition}

A natural question we can ask is how this theory relates to the theory
of Vassiliev invariants.

\begin{theorem} \label{s:th1}
Every Vassiliev invariant is a linking class finite type invariant
for all classes.
\end{theorem}

\begin{proof}
The theorem follows from the equality
\[ 
  \printname{figure15ab}
  \setlength{\unitlength}{0.5\standardunitlength}
  \begin{array}{c}  \hspace{-1.7mm}
    \raisebox{-8pt}{\input draws/figure15ab.tex }
    \hspace{-1.9mm}
  \end{array}
 \]
Hence the value of $V$ on a DD($\bigwedge$)-singular P-link having
$m+1$ pairs of double points descends to the sum of values of $V$
on singular links each of which has $m+1$ double points.
\end{proof}

\begin{theorem}
If $V$ is a Vassiliev
invariant of pure tangles which descends to a well defined
invariant of links within a given linking class, via the closure
map, then $V$ in an invariant of type $m$, in the sense presented
here, for the links in that class.
\end{theorem}

The proof of this theorem uses exactly the same argument as the
previous proof, but before we calculate the value of $V$ on some link we
first cut it open to a pure tangle and take the value of $V$ on
this tangle.

\begin{cor}
Since the Milnor $\mu_{ijk}$ is an invariant of pure tangles which is well
defined as an invariant of algebraically split links, and by~\cite{ba:2,
Lin:Expansions} is known to be an invariant of type $2$, we obtain as
a corollary that it is an invariant of at most type $2$ in our sense too.
\end{cor}

\begin{theorem} \label{howmany} (proof in Section~\ref{sec:low})
The only invariants of finite type for links up to equivalence class
in the sense presented here, up to type $2$, are restrictions of
Vassiliev invariants and $\mu_{ijk}$.
\end{theorem}

The paper is organized as follows. In Section~\ref{sec:diagrams}
we introduce the spaces $\mathcal{D}^{DD}$ of Double Dating
diagrams and $\mathcal{D}^{\bigwedge }$ of Wedge diagrams and we
show that these spaces are the appropriate analogous of chord
diagrams. In Section~\ref{sec:relations} we introduce some
relations which are satisfied by double dating and wedge diagrams.
In Section~\ref{sec:low} we introduce the spaces $\mathcal{F}_P^m$
of all type $m$ invariants, and $\mathcal{A}^{DD}$ which is the
quotient space of $\mathcal{D}^{DD}$ divided by the subspace
generated by all relations we present in
Section~\ref{sec:relations}.  Furthermore in
Section~\ref{sec:low}, we calculate the dimension of those spaces
up to degree $2$, with respect to the number of components.  In
Section~\ref{sec:dd2c} we introduce a map from $\mathcal{D}^{DD}$
to the algebra of chord diagram, from Vassiliev invariants, and we
describe the image under this map of $\mathcal{D}^{DD},
\mathcal{D}^{\bigwedge}$.
 Finally, we present some open questions.

At the time that this preprint was still in preparation I
received a preprint~\cite{me} by B.~Mellor, which is based on an
earlier version of my preprint. Thus there is some overlap
between the two works. The basic definitions are nearly the same,
and in both places it is shown that the Milnor invariant
$\mu_{ijk}$ is an invariant of finite type for algebraically split
links.
 The proofs are
different though.
\subsection{Acknowledgement}
The theory presented here is a joint work with Dror Bar-Natan from
the
Hebrew University in Jerusalem. I am deeply grateful to Dror for
his sincere help and advisory during the time this work was done
and all along
my Ph.D studies.

I also wish to express my deep gratitude to Stavros
Garoufalidis for his warm hospitality and the stimulating atmosphere
and useful conversations I had with him while visiting him at
Harvard University during December 98.

\section{The diagrams} \lbl{sec:diagrams}
We wish to present here the objects appropriate to play the role of chord
diagrams in each of the two types of iteration processes we have
described above. First, we will treat the  case of DD-adjacent links, in which
Double Dating diagrams (DD-diagrams, for short) will be presented.
Then we tend to the case of $\bigwedge$-adjacent links, and introduce Wedge diagrams
\subsection{Double Dating Diagrams}
\begin{figure}
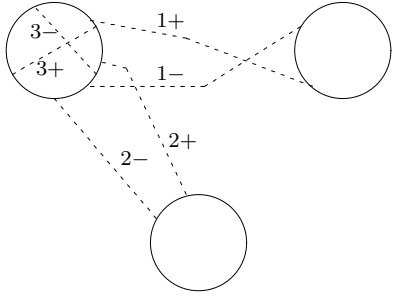

\[ 
  \printname{figure20}
  \setlength{\unitlength}{0.5\standardunitlength}
  \begin{array}{c}  \hspace{-1.7mm}
    \raisebox{-8pt}{\input draws/figure20.tex }
    \hspace{-1.9mm}
  \end{array}
 \]
\caption{A double dating diagram} \label{dd-diag}
\end{figure}

\begin{definition}
A degree $m$ double dating diagram on $l$ components is a collection of
$l$ ordered oriented circles, and $m$ ordered pairs, ordered by $+/-$
signs, of dashed lines(chords) so that both chords in each pair have
their end points on the same one/two circles. In Figure~\ref{dd-diag}
an example of a double dating diagram of degree $3$ on three components
is shown. Note that the numerals 1-3 written next to the chords, in
the figure, are not part of the structure of a DD-diagrams. They are
only a mean of presenting the pairing of the chords.
\end{definition}
Let $\mathcal{D}^{DD}_{m}$ be the span of all double dating diagrams of
degree $m$ and let $\mathcal{D}^{DD}$ be the direct sum over all $m$
of $\mathcal{D}^{DD}_{m}$.
\begin{theorem} \label{th2.1}
The space $\mathcal{D}^{DD}_{m}$ is analogous to the space of chord
diagrams of degree $m$, defined in the theory of Vassiliev invariants. To
be more specific,
\begin{enumerate}
\item \label{2.1.i}
To every degree $m$ DD-singular link, L, there corresponds a degree $m$
DD diagram, $D_{L}$.
\item
\begin{enumerate}
\item \label{2.1.ii}
To every degree $m$ DD diagram, D, there  exists a degree
 $m$ singular link, L(D), (in every linking class) so that
\begin{equation}
D_{L(D)} = D
\end{equation}
\item \label{2.1.iii}
Two such links, within the same class, corresponding to the same DD
diagram, differ by finite sequence of double crossing changes.
\end{enumerate}
\end{enumerate}
\end{theorem}

\begin{proof}

The proof of~\ref{2.1.i} is best illustrated by an example shown in
Figure~\ref{s:link}, where we show a singular link on the
left hand side and its corresponding double dating diagram on the right
hand side.  For the proof of~\ref{2.1.ii} we show how to build a singular
link for a given DD diagram. Let D be a DD diagram with some full lines
and pairs of chords. We will call the chords belonging to a certain pair
``pals''. In Figure~\ref{a-b} we see a double dating diagram on the left
hand side.  Use ``finger moves'' along the chords in order to get from (a)
to (b).  Following this procedure to all pairs of pals chords will give
a singular link as desired. One may think of the finger move as axon
connecting two neurons by a synapse. Clearly this link is not unique
and it can be obtained in any linking class. However now we show that
any two singular links respecting the same DD diagram and within the
same class differ by a finite sequence of double crossing flips. This
is what is claimed in~\ref{2.1.iii}.
\begin{figure}
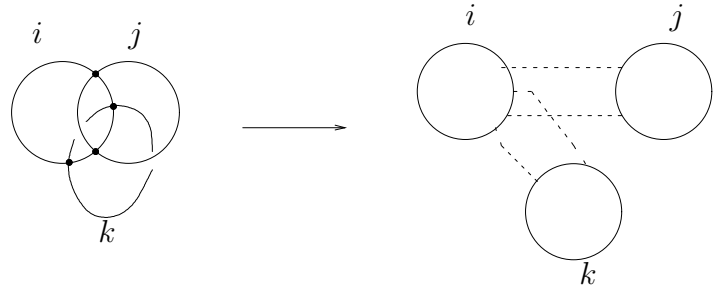

\[ 
  \printname{figure14}
  \setlength{\unitlength}{0.5\standardunitlength}
  \begin{array}{c}  \hspace{-1.7mm}
    \raisebox{-8pt}{\input draws/figure14.tex }
    \hspace{-1.9mm}
  \end{array}
 \]
\caption{A singular link and its DD-diagram} \label{s:link}
\end{figure}
\begin{figure}
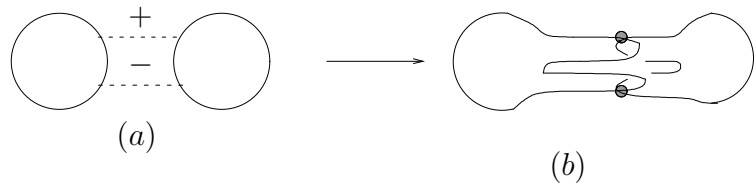

\[ 
  \printname{figure21b}
  \setlength{\unitlength}{0.5\standardunitlength}
  \begin{array}{c}  \hspace{-1.7mm}
    \raisebox{-8pt}{\input draws/figure21b.tex }
    \hspace{-1.9mm}
  \end{array}
 \]
\caption{From a diagram to a singular link} \label{a-b}
\end{figure}
We say that a link (singular or not), is in a neural position if
it is embedded in $S^{3}$ so that each component is an unknotted
circle remote from all other components, knotted and linked to
other components by a network of axons and synapses. Every link is
isotopic to a link in a neural position. To see this one can take
a disk neighborhood of every crossing of a given link diagram. The
crossing can be regular or singular. Making the disk thicker one
gets a $3$-ball neighborhood of the crossing. It is possible to
move this ball-neighborhood around in $3$-space thus creating an
axon of the neuro-link and the crossing neighborhood becomes the
synapse.
A self intersection point belonging to a singular
link will give a singular synapse . A non singular synapse, is either positive
or negative according to its contribution ($+1$ for positive, $-1$ for
negative) to the linking matrix.  See Figure~\ref{n-pos} for examples
of a link and a knot in neural positions, and Figure~\ref{synp} for the
distinguishing between the signs of non-singular synapse and also for an
example of a singular synapse.
\begin{figure}
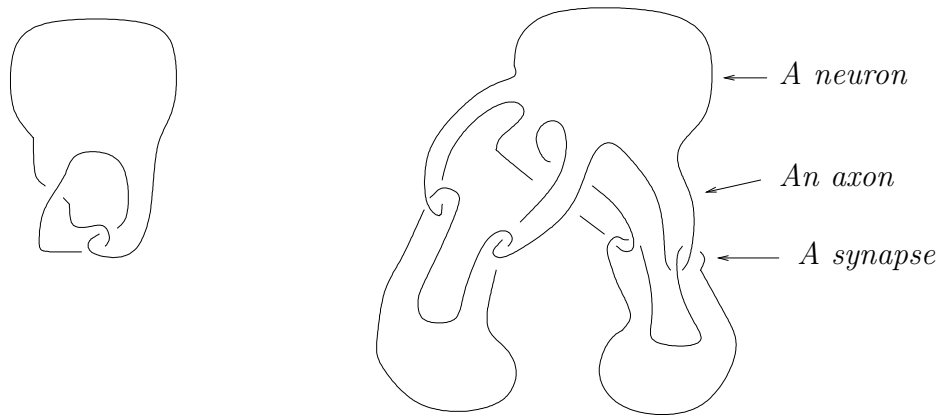

\[ 
  \printname{figure21c}
  \setlength{\unitlength}{0.5\standardunitlength}
  \begin{array}{c}  \hspace{-1.7mm}
    \raisebox{-8pt}{\input draws/figure21c.tex }
    \hspace{-1.9mm}
  \end{array}
 \]
\caption{A neuro trefoil knot and a neuro Borromean link} \label{n-pos}
\end{figure}
\begin{figure}
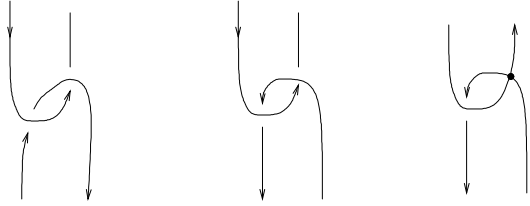

\[ 
  \printname{figure21ca}
  \setlength{\unitlength}{0.5\standardunitlength}
  \begin{array}{c}  \hspace{-1.7mm}
    \raisebox{-8pt}{\input draws/figure21ca.tex }
    \hspace{-1.9mm}
  \end{array}
 \]
\caption{
  A positive synapse, a negative synapse and a singular synapse (in order)
} \label{synp}
\end{figure}
Given two links, $L_1, L_2$, respecting
the same DD diagram, D and the same linking matrix, we first bring both
of them to neural positions.
Suppose $i_1,i_2$ $j_1,j_2$ are the $i^{th}$ and $j^{th}$ components
of $L_1, L_2$
respectively. The numbers $i$ and $j$ of both links are given by the diagram
$D$.
 First, we will make sure that modulo double crossing change,
both $i_1,j_1$ and $i_2,j_2$ have the same
number of synapse involving each pair.
 Since the links respect the same diagram we
do not have to worry about singular synapses. The number of them is given from
the diagram. suppose that the pair $i_2,j_2$ has more non-singular synapses than
$i_1,j_1$. Since the links $L_1, L_2$ have the same linking matrix, and the fact
that each synapse adds $\pm 1$ to the linking matrix we have that,
\begin{displaymath}
\sharp \{\textrm{synapses between $i_2,j_2$}\}  -\sharp \{\textrm{synapses
between $i_1,j_1$}\} = \textrm{ an even number}
\end{displaymath}
It follows that we can find a set of synapses involving $i_2,j_2$
representing this difference which can be divided into pairs, each
of them consists of one positive and one negative synapse. Every
such pair gives us a double crossing flip that allows us to cancel
those synapses and hence decreases the number of synapses between
$i_2,j_2$ by $2$. We repeat this argument for every pair of
components of $L_2$. This way we change $L_2$ into a link,
$L_{2}'$, every two components of which has the same number of
synapses in the corresponding pair of components in $L_1$. In
order to finish the proof we still have to show that we can find a
sequence of double flips making $L_{2}'$ The same as $L_1$. Since
an axon is made of two strings going in two opposite directions it
is transparent to other axons and circles in the sense that we can
flip it through axons and circles using only double crossing
flips. Therefore the combinatorics of the setting which is defined
entirely by the diagram determines the topology.
\[
  \begin{array}{ccc}
    
  \printname{figure21d1}
  \setlength{\unitlength}{0.5\standardunitlength}
  \begin{array}{c}  \hspace{-1.7mm}
    \raisebox{-8pt}{\input draws/figure21d1.tex }
    \hspace{-1.9mm}
  \end{array}
 & \qquad & 
  \printname{figure21d2}
  \setlength{\unitlength}{0.5\standardunitlength}
  \begin{array}{c}  \hspace{-1.7mm}
    \raisebox{-8pt}{\input draws/figure21d2.tex }
    \hspace{-1.9mm}
  \end{array}
 \\
    \ && \\
    \parbox{2.6in}{an axon of the $i^{th}$ component over-crosses an axon
      of the $j^{th}$ component
    } & &
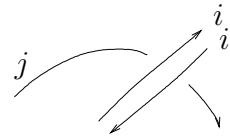
    \parbox{2.6in}{an axon of the $i^{th}$ component over-crosses the
      $j^{th}$ circle component
    }
  \end{array}
\]
The links $L_1$ and $L_2'$ have the same combinatorial setting and
the axons of $L_1$ can be moved freely, (up to double crossing
flips), in space onto the axons of $L_2'$ so they become isotopic
links. This completes the proof.
\end{proof}
\subsection{Wedge Diagrams}
\begin{figure}
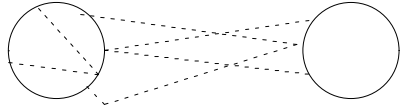

\[ 
  \printname{figure15}
  \setlength{\unitlength}{0.5\standardunitlength}
  \begin{array}{c}  \hspace{-1.7mm}
    \raisebox{-8pt}{\input draws/figure15.tex }
    \hspace{-1.9mm}
  \end{array}
 \]
\caption{A $\bigwedge$-diagram} \label{s:w-diag}
\end{figure}

\begin{definition} A wedge diagram of degree $m$ on $l$ components
is a collection of $l$ oriented ordered circles and $m$ dashed wedges, so that
the tip of a wedge lies on one circle and both its legs lies on another
(maybe the same) circle. The legs of each wedge are ordered by $+/-$ signs.
See Figure~\ref{s:w-diag}, for an example of a
Wedge diagram of degree $3$ and $2$ components.
\end{definition}
Define $\mathcal{D}^{\bigwedge}_{m}$ and $\mathcal{D}^{\bigwedge}$
in the same way as $\mathcal{D}^{DD}$ and $\mathcal{D}^{DD}_{m}$
\begin{theorem}
The space of chord diagrams corresponding to $\bigwedge$-singular links
(within a linking class), in the sense of theorem~\ref{th2.1},
is the space $\mathcal{D}^{\bigwedge}$.
\end{theorem}
\begin{proof}
The proof of this theorem is basically identical to that
of~\ref{th2.1}. The only difference occurs when we reduce the number of
axons in $L_{2}$ so that it would be the same as for $L_{1}$.  as already
mentioned there, we cancel a pair of axons by the use of double crossing
flips.  In general these crossings do not have to be close. If the two
axons are adjacent then those crossings are close. If they are not and
we have two axons separated by a third one, as below, then we slide one of
them as indicated and make the double crossing flip encircled which is
a wedge type flip, so that they become adjacent.
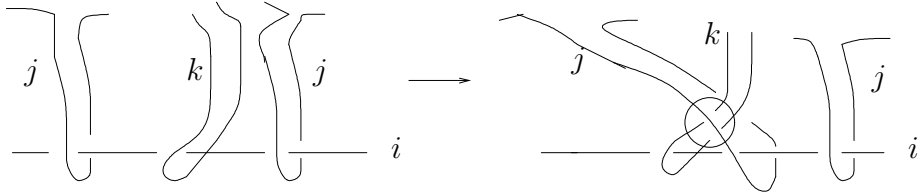
\[ 
  \printname{figure21e}
  \setlength{\unitlength}{0.5\standardunitlength}
  \begin{array}{c}  \hspace{-1.7mm}
    \raisebox{-8pt}{\input draws/figure21e.tex }
    \hspace{-1.9mm}
  \end{array}
 \]
\end{proof}
Wedge diagrams naturally mapped into Double Dating diagrams as we
can put the end points at the tip of each wedge further from each
other. This is actually the illustration of the inclusion map from
$\mathcal{D}^{\bigwedge}$ into $\mathcal{D}^{DD}$. At the end of
the next section we will show that this map is onto.

\section{Relations For Diagrams} \lbl{sec:relations}

In order to define weight systems in this theory, as well as building
a Hutching's theory for giving integration conditions of weight systems
to invariants~\cite{hu:1}, we need to have complete knowledge of the
relations satisfied by double dating and wedge diagrams. In this section
we will present the relations we already know. We suspect that the list
of relations below is complete, but we have not proven that yet. Recall
from~\cite{ba:1} that a weight system is a linear functional on the
space freely generated by chord diagrams which corresponds to finite type
invariants of links. In order to be a weight system oppose to a general
linear functional it is necessary that the functional will satisfy certain
relations originated by the topological nature of the setting. In the
following description of relations we present the topological origin and
translate it into relations among DD-diagrams ($\bigwedge$-diagrams). A
weight system then will become a linear functional on diagrams which
vanishes on the subspace spanned by all the relations.

\subsection{DD Diagrams}
\subsubsection{2-T Relation} In the following picture we can move the loop
circling the double point upward, using double crossing flips.
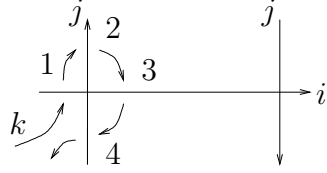
\[ 
  \printname{figure220}
  \setlength{\unitlength}{0.5\standardunitlength}
  \begin{array}{c}  \hspace{-1.7mm}
    \raisebox{-8pt}{\input draws/figure220.tex }
    \hspace{-1.9mm}
  \end{array}
 \]
Change crossings $1$ and $3$ simultaneously and then $2$ and $4$.
The result now will be the 2-T relation.
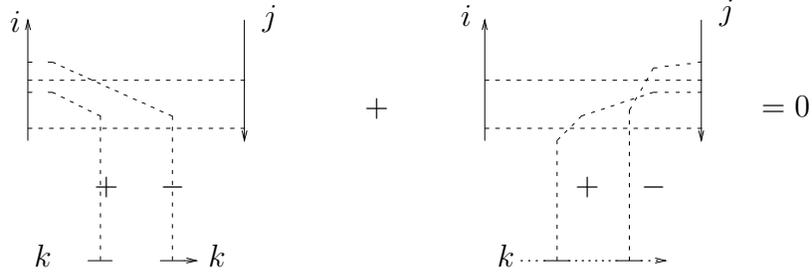
\[ 
  \printname{figure22a}
  \setlength{\unitlength}{0.5\standardunitlength}
  \begin{array}{c}  \hspace{-1.7mm}
    \raisebox{-8pt}{\input draws/figure22a.tex }
    \hspace{-1.9mm}
  \end{array}
 \]
\subsubsection{3-T Relation} Consider the following sequence of moves
involving the $i^{th}$ and $j^{th}$ components of the link.
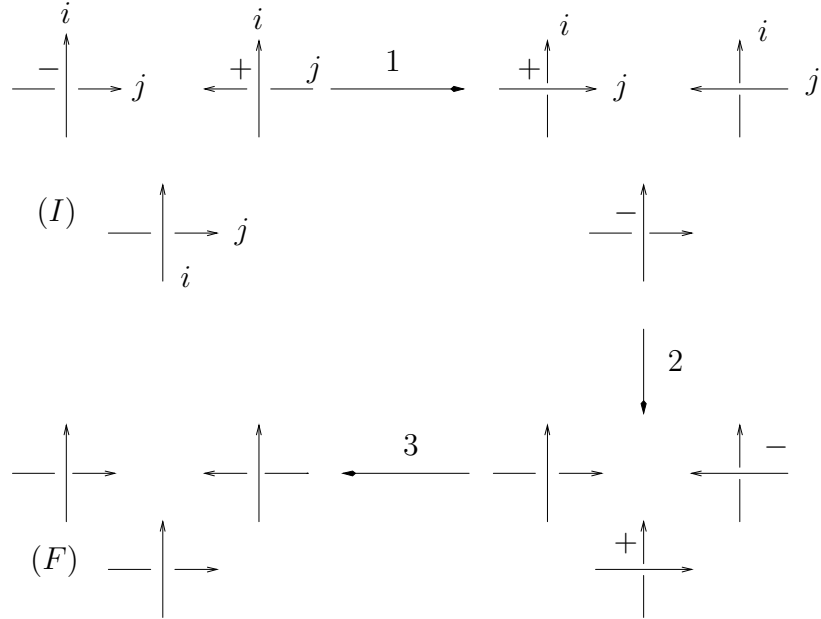
\[ 
  \printname{figure24}
  \setlength{\unitlength}{0.5\standardunitlength}
  \begin{array}{c}  \hspace{-1.7mm}
    \raisebox{-8pt}{\input draws/figure24.tex }
    \hspace{-1.9mm}
  \end{array}
 \]
As a result of applying this sequence, the initial and final
configurations (I) and (F) are the same. Then in terms of the diagrams
this sequence translates into the fact that each pair of chords created while
performing each step can be rewritten in terms of the other two
pairs.
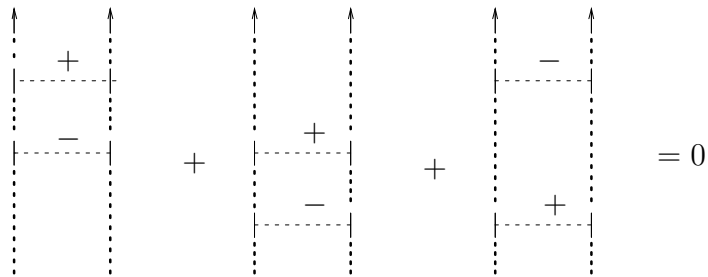
\[ 
  \printname{figure24a}
  \setlength{\unitlength}{0.5\standardunitlength}
  \begin{array}{c}  \hspace{-1.7mm}
    \raisebox{-8pt}{\input draws/figure24a.tex }
    \hspace{-1.9mm}
  \end{array}
 \]
As an immediate consequence of the 3-T relation we obtain a
1-T relation, which also can be proven using
the second Reidemeister move.
\[ 
  \printname{figure27a}
  \setlength{\unitlength}{0.5\standardunitlength}
  \begin{array}{c}  \hspace{-1.7mm}
    \raisebox{-8pt}{%
\begingroup\makeatletter\ifx\SetFigFont\undefined
\def\x#1#2#3#4#5#6#7\relax{\def\x{#1#2#3#4#5#6}}%
\expandafter\x\fmtname xxxxxx\relax \def\y{splain}%
\ifx\x\y   
\gdef\SetFigFont#1#2#3{%
  \ifnum #1<17\tiny\else \ifnum #1<20\small\else
  \ifnum #1<24\normalsize\else \ifnum #1<29\large\else
  \ifnum #1<34\Large\else \ifnum #1<41\LARGE\else
     \huge\fi\fi\fi\fi\fi\fi
  \csname #3\endcsname}%
\else
\gdef\SetFigFont#1#2#3{\begingroup
  \count@#1\relax \ifnum 25<\count@\count@25\fi
  \def\x{\endgroup\@setsize\SetFigFont{#2pt}}%
  \expandafter\x
    \csname \romannumeral\the\count@ pt\expandafter\endcsname
    \csname @\romannumeral\the\count@ pt\endcsname
  \csname #3\endcsname}%
\fi
\fi\endgroup
{\renewcommand{\dashlinestretch}{30}
\begin{picture}(2554,1839)(0,-10)
\path(12,1062)(12,462)
\path(1812,1062)(1812,462)
\dottedline{90}(12,1137)(12,1812)
\path(42.000,1692.000)(12.000,1812.000)(-18.000,1692.000)
\dottedline{90}(1812,1137)(1812,1812)
\path(1842.000,1692.000)(1812.000,1812.000)(1782.000,1692.000)
\dottedline{90}(12,387)(12,12)
\dottedline{90}(1812,387)(1812,12)
\put(2262,612){\makebox(0,0)[lb]{\smash{{\mathmode{= 0 }}}}}
\dashline{60.000}(12,687)(1812,687)
\dashline{60.000}(12,912)(1812,912)
\end{picture}
}
 }
    \hspace{-1.9mm}
  \end{array}
 \]

\subsection{ $\bigwedge$ Diagrams}
First note that the 2-T, 3-T and their consequences presented earlier
are also valid for wedge diagrams.  the following two figures show the
appropriate setting of those relations for wedge diagrams.

\subsubsection{$\bigwedge$-2-T}
for wedge diagrams the 2-T relation becomes,
\[ 
  \printname{figure27}
  \setlength{\unitlength}{0.5\standardunitlength}
  \begin{array}{c}  \hspace{-1.7mm}
    \raisebox{-8pt}{\input draws/figure27.tex }
    \hspace{-1.9mm}
  \end{array}
 \]
Notice the change of components where the wedge type intersection
occurs.

\subsubsection{$\bigwedge$-3-t}
\[ 
  \printname{figure27b}
  \setlength{\unitlength}{0.5\standardunitlength}
  \begin{array}{c}  \hspace{-1.7mm}
    \raisebox{-8pt}{\input draws/figure27b.tex }
    \hspace{-1.9mm}
  \end{array}
 \]

\subsubsection{Interchange Relation} The situation described in the
next figure can be treated as in the 2-T relation, only that
here we must assure that the crossings we change are close.
\[ 
  \printname{figure90}
  \setlength{\unitlength}{0.5\standardunitlength}
  \begin{array}{c}  \hspace{-1.7mm}
    \raisebox{-8pt}{\input draws/figure90.tex }
    \hspace{-1.9mm}
  \end{array}
 \]
This
can be done in two inequivalent ways, each corresponds to a
certain pairing of the crossings.
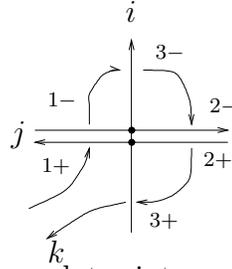
\begin{enumerate}
\item \hfill$
  \printname{figure25ab}
  \setlength{\unitlength}{0.5\standardunitlength}
  \begin{array}{c}  \hspace{-1.7mm}
    \raisebox{-8pt}{\input draws/figure25ab.tex }
    \hspace{-1.9mm}
  \end{array}
$\hfill\ \newline
In terms of wedge diagrams, this translates into an interchange relation.
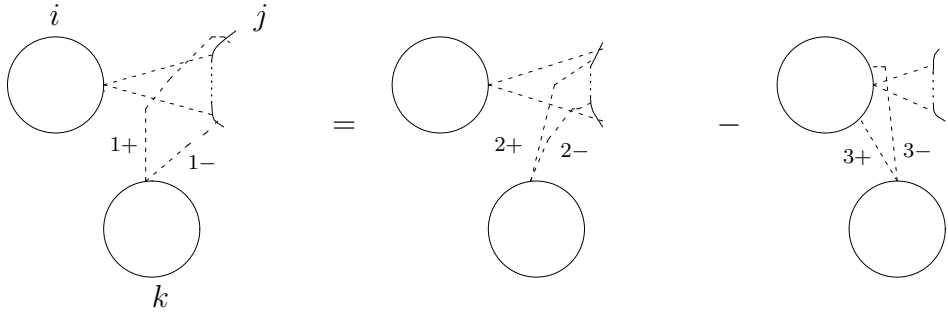
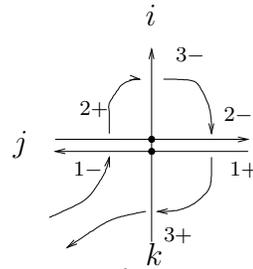
\[ 
  \printname{figure25ac}
  \setlength{\unitlength}{0.5\standardunitlength}
  \begin{array}{c}  \hspace{-1.7mm}
    \raisebox{-8pt}{\input draws/figure25ac.tex }
    \hspace{-1.9mm}
  \end{array}
 \]
\item \hfill$
  \printname{figure25}
  \setlength{\unitlength}{0.5\standardunitlength}
  \begin{array}{c}  \hspace{-1.7mm}
    \raisebox{-8pt}{\input draws/figure25.tex }
    \hspace{-1.9mm}
  \end{array}
$\hfill\ \newline
In this case we get the following version of the interchange relation.
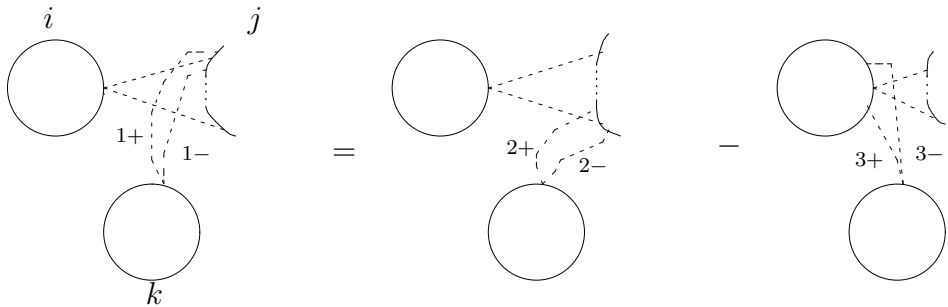
\hfill
\[ 
  \printname{figure25a}
  \setlength{\unitlength}{0.5\standardunitlength}
  \begin{array}{c}  \hspace{-1.7mm}
    \raisebox{-8pt}{\input draws/figure25a.tex }
    \hspace{-1.9mm}
  \end{array}
 \] 
\end{enumerate}

\par\noindent\parbox{4.2in}{
  \subsubsection{Wood Pecker relation} Now consider a case where the
  two pieces of the $k^{th}$ component loop around a double singularity
  as shown on the right.  Again move the loops upward, performing the
  necessary double crossing changes and follow them by the matching
  wedge diagrams. You will get the Wood Pecker relation shown below:
}\hfill$
  \printname{figure26}
  \setlength{\unitlength}{0.7\standardunitlength}
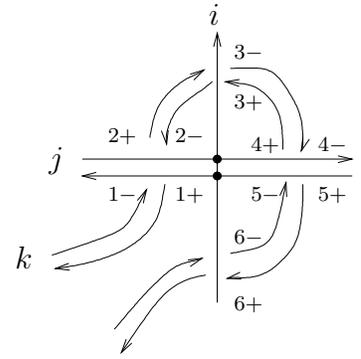
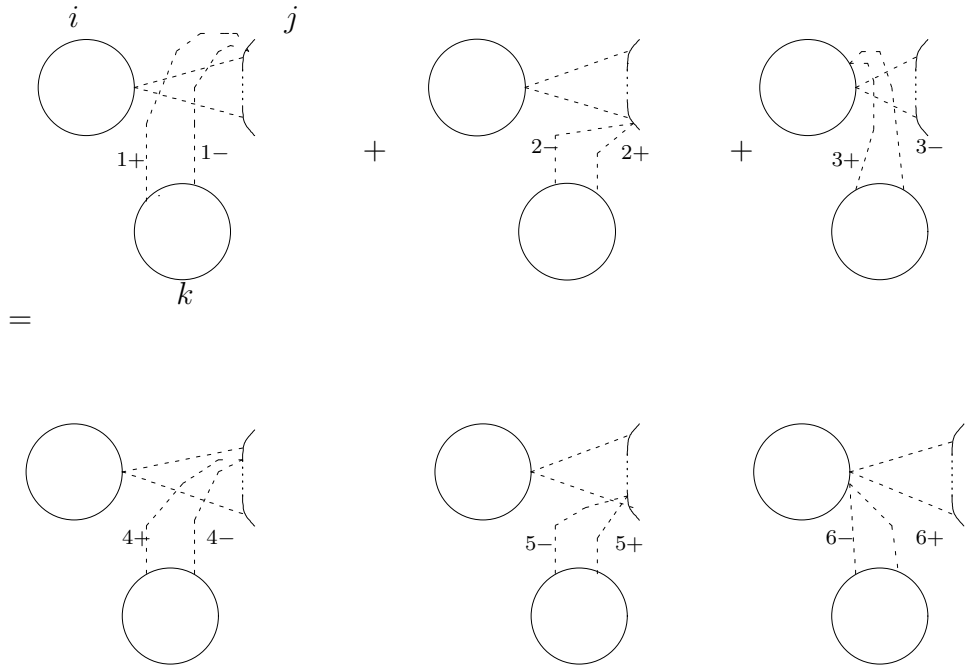
  \begin{array}{c}  \hspace{-1.7mm}
    \raisebox{-8pt}{\input draws/figure26.tex }
    \hspace{-1.9mm}
  \end{array}
$\hfill
\[ 
  \printname{figure26a}
  \setlength{\unitlength}{0.5\standardunitlength}
  \begin{array}{c}  \hspace{-1.7mm}
    \raisebox{-8pt}{\input draws/figure26a.tex }
    \hspace{-1.9mm}
  \end{array}
 \]
As said earlier in section $2$ 
Wedge diagrams naturally included in Double Dating diagrams. Since
diagrams are the dual space of the space of invariants and since
every DD-finite type invariant is a $\bigwedge$-finite type
invariant this  map between the spaces of diagrams should be onto.
Indeed it is. The topological proof of this fact is illustrated in
the following.
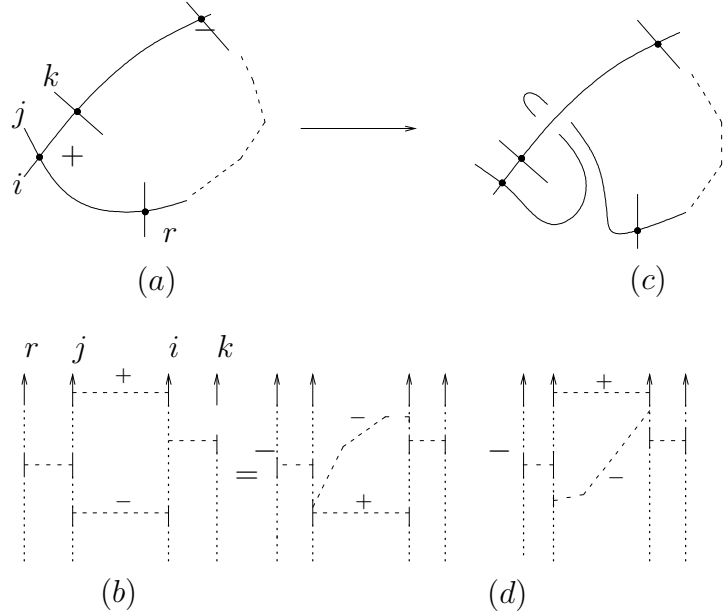
\[ 
  \printname{figure01}
  \setlength{\unitlength}{0.5\standardunitlength}
  \begin{array}{c}  \hspace{-1.7mm}
    \raisebox{-8pt}{\input draws/figure01.tex }
    \hspace{-1.9mm}
  \end{array}
 \]
In the above left (a) we have a pair of intersection points
 between the components $i$ and$j$ signed $+/-$,and divided by intersections
with the $k$ and $r$ components.
 (b) is the double dating diagram corresponding to it.
We use a finger move (c) and the $3$-T relation in order to have (b),
as a sum of diagrams (d) in which one of them is a wedge diagram and the other
is ``almost'' wedge. Applying this procedure as much as necessary makes the
double dating diagram (b) become a sum of wedge diagrams.

\section{ Some Low Degree Computations} \lbl{sec:low}

Denote the vector space of all type $n$ invariants of links with $m$ components,
and linking matrix $P$, by $\mathcal{F}^n_{P,m}$.
\begin{theorem}
The dimensions of $\mathcal{F}^m_{P,m}$ for $n = \{ 0,1,2\}$
are given by,\\
dim $\mathcal{F}^0_{P,m} = 1$ for all $m$.\\
dim $\mathcal{F}^1_{P,m} = 0$ for all $m$.\\
dim $\mathcal{F}^2_{P,m} = m$ for all $m$ where $P \neq 0$.\\
dim $\mathcal{F}^2_{0,3} = 3 + 1 $.\\
dim $\mathcal{F}^2_{0,m} = m +\binom{m}{3}$ for all $m>3$.
\end{theorem}
\begin{proof}
The invariants of type $0$ are the constants.
 There is a unique Vassiliev invariant of type $2$ for knots, which is
$V_2$ (the second coefficient of the Conway invariant). Therefore
for a link with $m$ component, we have $m$ $V_2's$ each one for
each component of the link. The extra invariant in
$\mathcal{F}^2_{0,3}$ is the Milnor $\mu_{ijk}$. If $m>3$ then we
also have $\binom{m}{3}$ such invariants. The list of finite type
invariants given above, consists of the invariant of type $\leq 2$
which come from Vassiliev invariant of links and tangles.
This list gives lower bounds on the dimensions claimed in the
theorem.  The completion of the proof is done by calculating the
dimensions of the space of double dating diagrams up to degree
$2$. This calculation is the content of the next theorem.
\end{proof}
Let $\mathcal{A}^{DD,n}_m$ be $\mathcal{D}^{DD,n}_m$ modulo
1-T, 2-T, 3-T.
\begin{theorem} \label{s:lowdim}
The dimensions of $\mathcal{A}^{DD,n}_m$ up to degree 2 are given by,\\
dim $\mathcal{A}^{DD,0}_m = 1$ for all $m$.\\
dim $\mathcal{A}^{DD,1}_m = 0$ for all $m$.\\
dim $\mathcal{A}^{DD,2}_m = \binom{m}{3}$
\end{theorem}

\begin{proof}
The proof of the theorem is done by using the
1-T, 2-T and 3-T  which give sequences of equalities
between double dating diagrams. Those equalities are given in
Appendix A.
\end{proof}

\section{From Double Dating Diagrams to Chord Diagrams} \lbl{sec:dd2c}
We briefly recall the definitions of the algebra
$\mathcal{A}^{c}$ of chord diagrams, and
$\mathcal{A}^{t}$ of trivalent diagrams. See~\cite{ba:1} for
further details. A chord diagram is a disjoint union of oriented
circles with dashed chords attached to them.
\begin{displaymath}
\mathcal{A}^{c} = \textrm{span} \{ \textrm{chord diagrams} \} /
\textrm{4-T relation}
\end{displaymath}
\begin{figure}
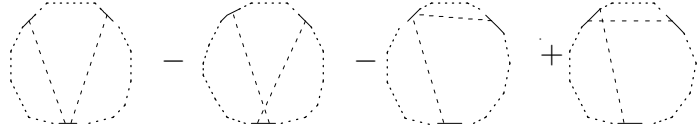

\[ 
  \printname{figure70a}
  \setlength{\unitlength}{0.5\standardunitlength}
  \begin{array}{c}  \hspace{-1.7mm}
    \raisebox{-8pt}{\input draws/figure70a.tex }
    \hspace{-1.9mm}
  \end{array}
 \]
\caption{4-T relation}
\end{figure}
A trivalent diagram is the same as chord diagram with the
exception that there are internal vertices allowed, where an
internal vertex is a vertex where three chords meet. Every
internal vertex is oriented by ordering the three edges meeting at
this vertex. $\mathcal{A}^{t}$ is the algebra generated by
trivalent diagrams modulo the previous 4-T, plus the additional AS
and STU relations.
\[ 
  \printname{figure70b}
  \setlength{\unitlength}{0.5\standardunitlength}
  \begin{array}{c}  \hspace{-1.7mm}
    \raisebox{-8pt}{\input draws/figure70b.tex }
    \hspace{-1.9mm}
  \end{array}
 \]
In~\cite{ba:1} it is shown that $\mathcal{A}^c$ and $\mathcal{A}^t$ are
isomorphic.
Recall from Theorem ~\ref{s:th1} that Vassiliev invariants are naturally mapped
into finite type invariants under the inclusion map. A
natural question arising from this observation is whether or not
this map is an $1:1$ and is it onto.
In this section we will introduce a map from the space
$\mathcal{A}^{DD}$ into $\mathcal{A}^{c}$ which is the natural dual map to the
inclusion map from Vassiliev invariants into finite type invariants. We shell
find the image under this map of the space $\mathcal{A}^{\bigwedge}$ of wedge
diagrams and since the inclusion map from wedge diagram to double dating diagram
was shown to be onto, we can identify also the image of $\mathcal{A}^{DD}$.
\begin{definition}
Define $\iota : \mathcal{D}^{DD} \to \mathcal{D}^{c}$ by defining
\[ 
  \printname{figure31}
  \setlength{\unitlength}{0.5\standardunitlength}
  \begin{array}{c}  \hspace{-1.7mm}
    \raisebox{-8pt}{\input draws/figure31.tex }
    \hspace{-1.9mm}
  \end{array}
 \]
On wedge diagrams the restriction of $\iota$ will look as
\[ 
  \printname{figure31a}
  \setlength{\unitlength}{0.5\standardunitlength}
  \begin{array}{c}  \hspace{-1.7mm}
    \raisebox{-8pt}{%
\begingroup\makeatletter\ifx\SetFigFont\undefined
\def\x#1#2#3#4#5#6#7\relax{\def\x{#1#2#3#4#5#6}}%
\expandafter\x\fmtname xxxxxx\relax \def\y{splain}%
\ifx\x\y   
\gdef\SetFigFont#1#2#3{%
  \ifnum #1<17\tiny\else \ifnum #1<20\small\else
  \ifnum #1<24\normalsize\else \ifnum #1<29\large\else
  \ifnum #1<34\Large\else \ifnum #1<41\LARGE\else
     \huge\fi\fi\fi\fi\fi\fi
  \csname #3\endcsname}%
\else
\gdef\SetFigFont#1#2#3{\begingroup
  \count@#1\relax \ifnum 25<\count@\count@25\fi
  \def\x{\endgroup\@setsize\SetFigFont{#2pt}}%
  \expandafter\x
    \csname \romannumeral\the\count@ pt\expandafter\endcsname
    \csname @\romannumeral\the\count@ pt\endcsname
  \csname #3\endcsname}%
\fi
\fi\endgroup
{\renewcommand{\dashlinestretch}{30}
\begin{picture}(7133,3481)(0,-10)
\put(5175,1808){\makebox(0,0)[lb]{\smash{{\mathmode{-}}}}}
\put(7162.500,1733.000){\arc{13425.838}{2.8817}{3.4015}}
\put(-4387.500,1733.000){\arc{13425.838}{6.0233}{6.5431}}
\put(1350,683){\ellipse{1200}{1200}}
\put(3825,683){\ellipse{1200}{1200}}
\put(6525,683){\ellipse{1200}{1200}}
\dashline{60.000}(975,2333)(1350,1283)
\dashline{60.000}(3525,2258)(3825,1283)
\dashline{60.000}(6900,2408)(6525,1283)
\dashline{60.000}(1350,1283)(1725,2333)
\dottedline{90}(1125,2183)(1350,2108)(1650,2258)
\dottedline{90}(3600,2258)(3825,2108)
\dottedline{90}(3450,2408)(3375,2633)
\dottedline{90}(6750,2258)(6600,2183)(6375,2183)
\dottedline{90}(6975,2558)(7050,2708)(7050,2783)
\dottedline{90}(1800,2408)(1875,2633)
\dottedline{90}(900,2408)(825,2558)(825,2708)
\path(1650,2258)(1800,2408)
\path(900,2408)(1125,2183)
\path(3450,2408)(3600,2258)
\path(6750,2258)(6975,2558)
\put(0,1808){\makebox(0,0)[lb]{\smash{{\mathmode{\iota}}}}}
\put(2625,1808){\makebox(0,0)[lb]{\smash{{\mathmode{=}}}}}
\put(900,1658){\makebox(0,0)[lb]{\smash{{\mathmode{\scriptstyle +}}}}}
\put(1650,1658){\makebox(0,0)[lb]{\smash{{\mathmode{\scriptstyle -}}}}}
\end{picture}
}
 }
    \hspace{-1.9mm}
  \end{array}
 \]
In general, $\iota$ of a diagram is defined by taking the
sum over all pairs of chords.
\end{definition}
It is easily checked that $\iota$ is well
defined map from $\mathcal{A}^{DD}$ to $\mathcal{A}^{c}$.
 This is done by showing that all the relations in
$\mathcal{A}^{DD}$ are mapped into the subspace of $\mathcal{A}^c$
spanned by all $4-T$ relations. This is done in appendix $B$.
We view $\iota$ as a map from $\mathcal{A}^{DD}$
into $\mathcal{A}^{t}$
\begin{theorem}
The  image of
$\mathcal{A}^{\bigwedge}$ under $\iota$ is the subspace, $\mathcal{A}^{st}$,of
$\mathcal{A}^t$, consisting of strutless diagrams in which there is
no chord having both its ends as external vertices.
\end{theorem}
\begin{proof}
In order to prove the theorem we introduce another space
$\mathcal{D}^{\bigwedge ,t}$ having $\mathcal{D}^{\bigwedge}$ as a
subspace. Further we will define a map $\overline \iota$ from
$\mathcal{D}^{\bigwedge ,t}$ into $\mathcal{A}^t$ which extends
$\iota$ and we will show that $\mathcal{D}^{\bigwedge ,t}$ is
mapped into the subspace of strutless diagrams. The objects of
$\mathcal{D}^{\bigwedge ,t}$ are diagrams having trivalent
vertices some of which external where a chord meets a circle
component of the diagram, and some are internal where three chords
meet. Some of the internal vertices are marked by $*$. The chords
that meet at an internal marked vertex are required that at least two of
them are attached to the same circle-component of the diagram.
those two chords are also signed by $+,-$ signs.
In Figure ~\ref{w,t diag} there is an example of a diagram
in $\mathcal{D}^{\bigwedge,t}$.
\begin{figure}
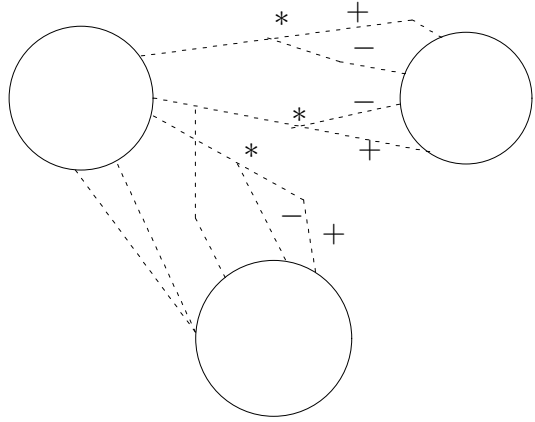

\[ 
  \printname{figure32}
  \setlength{\unitlength}{0.5\standardunitlength}
  \begin{array}{c}  \hspace{-1.7mm}
    \raisebox{-8pt}{\input draws/figure32.tex }
    \hspace{-1.9mm}
  \end{array}
 \]
\caption{A diagram in $\mathcal{D}^{\bigwedge,t}$}
\label{w,t diag}
\end{figure}
Now define $\overline \iota$ by
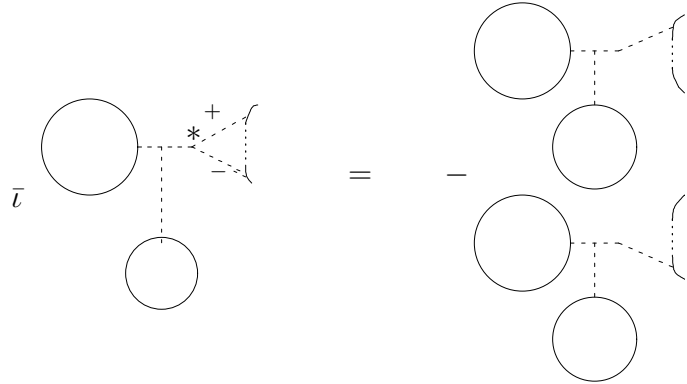
\[ 
  \printname{figure33}
  \setlength{\unitlength}{0.5\standardunitlength}
  \begin{array}{c}  \hspace{-1.7mm}
    \raisebox{-8pt}{\input draws/figure33.tex }
    \hspace{-1.9mm}
  \end{array}
 \]
and if a diagram contains trivalent vertices non of them is a marked vertex,
then it is mapped
to zero.
Notice that the following holds.
\[ 
  \printname{figure34}
  \setlength{\unitlength}{0.5\standardunitlength}
  \begin{array}{c}  \hspace{-1.7mm}
    \raisebox{-8pt}{\input draws/figure34.tex }
    \hspace{-1.9mm}
  \end{array}
 \]
Suppose now that there is a marked vertex for which there are
other strands touching the circle where the two end points are, on
the arc between those two points. Then after applying $\overline \iota$
as much as necessary, the following equality.
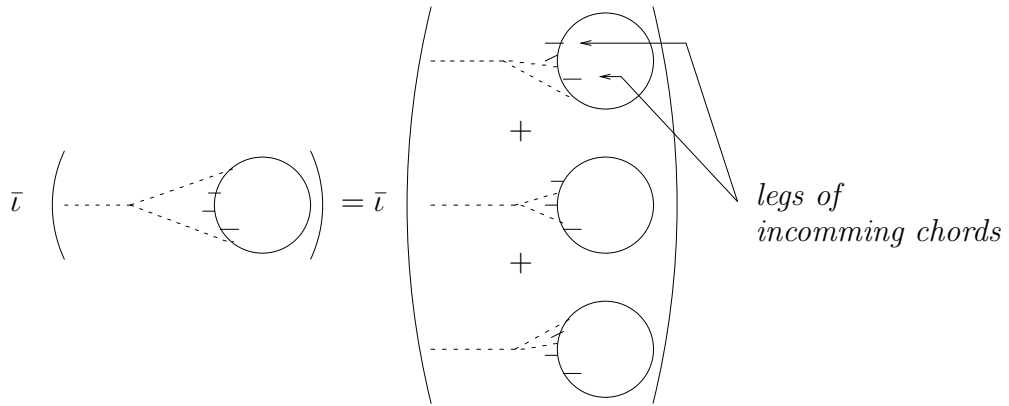
\[ 
  \printname{figure35}
  \setlength{\unitlength}{0.5\standardunitlength}
  \begin{array}{c}  \hspace{-1.7mm}
    \raisebox{-8pt}{\input draws/figure35.tex }
    \hspace{-1.9mm}
  \end{array}
 \]
Hence the only thing to worry about is an object as in the right
hand side of the previous equation. The claim is finally proven
by the two following equalities.
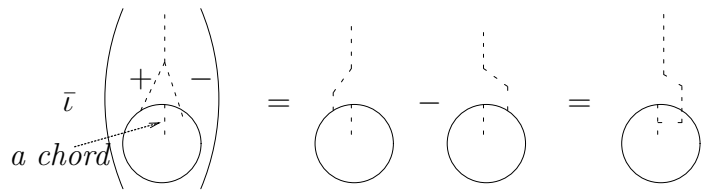
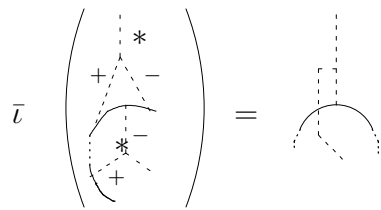
\[ 
  \printname{figure36}
  \setlength{\unitlength}{0.5\standardunitlength}
  \begin{array}{c}  \hspace{-1.7mm}
    \raisebox{-8pt}{\input draws/figure36.tex }
    \hspace{-1.9mm}
  \end{array}
 \]
\[ 
  \printname{figure37}
  \setlength{\unitlength}{0.5\standardunitlength}
  \begin{array}{c}  \hspace{-1.7mm}
    \raisebox{-8pt}{\input draws/figure37.tex }
    \hspace{-1.9mm}
  \end{array}
 \]
On the other hand, the fact that $\overline \iota$ is onto $\mathcal{A}^{st}$
 follows from the following equality.
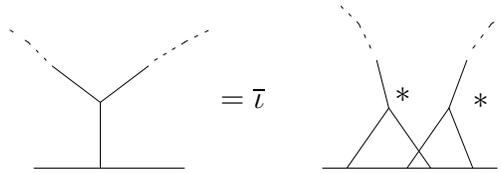
\[ 
  \printname{figure02}
  \setlength{\unitlength}{0.5\standardunitlength}
  \begin{array}{c}  \hspace{-1.7mm}
    \raisebox{-8pt}{%
\begingroup\makeatletter\ifx\SetFigFont\undefined
\def\x#1#2#3#4#5#6#7\relax{\def\x{#1#2#3#4#5#6}}%
\expandafter\x\fmtname xxxxxx\relax \def\y{splain}%
\ifx\x\y   
\gdef\SetFigFont#1#2#3{%
  \ifnum #1<17\tiny\else \ifnum #1<20\small\else
  \ifnum #1<24\normalsize\else \ifnum #1<29\large\else
  \ifnum #1<34\Large\else \ifnum #1<41\LARGE\else
     \huge\fi\fi\fi\fi\fi\fi
  \csname #3\endcsname}%
\else
\gdef\SetFigFont#1#2#3{\begingroup
  \count@#1\relax \ifnum 25<\count@\count@25\fi
  \def\x{\endgroup\@setsize\SetFigFont{#2pt}}%
  \expandafter\x
    \csname \romannumeral\the\count@ pt\expandafter\endcsname
    \csname @\romannumeral\the\count@ pt\endcsname
  \csname #3\endcsname}%
\fi
\fi\endgroup
{\renewcommand{\dashlinestretch}{30}
\begin{picture}(6249,2064)(0,-10)
\path(2187,12)(312,12)
\path(1137,12)(1137,837)
\path(1137,837)(537,1287)
\path(1137,837)(1737,1287)
\path(3912,12)(6087,12)
\path(4212,12)(4737,762)(5262,12)
\path(4962,12)(5487,762)(5787,12)
\path(4737,762)(4587,1362)
\path(5487,762)(5712,1362)
\dashline{60.000}(462,1362)(237,1587)(12,1737)
\dashline{60.000}(1812,1362)(2187,1587)(2487,1737)
\dashline{60.000}(4512,1512)(4362,1812)(4137,2037)
\dashline{60.000}(5787,1512)(6012,1812)(6237,1962)
\put(2637,762){\makebox(0,0)[lb]{\smash{{\mathmode{=  \overline \iota}}}}}
\put(4812,837){\makebox(0,0)[lb]{\smash{{\mathmode{*}}}}}
\put(5787,762){\makebox(0,0)[lb]{\smash{{\mathmode{*}}}}}
\end{picture}
}
 }
    \hspace{-1.9mm}
  \end{array}
 \]
Using STU we have that every strutless diagram is a sum of
diagrams in which every trivalent vertex is adjacent to a circle
as in the figure, and each such a diagram is in the image of
$\iota$
\end{proof}
\subsection{ The Case of Knots}
Let us focus on invariants of framed knots with a fixed self linking
number. Both double dating and chord diagrams in this case are made of
a single circle with dashed chords attached, and we have the
appropriate relations in each type of diagrams. Define a map $\nu$
from $\mathcal{A}^{c}$ into $\mathcal{A}^{DD}$ by setting,
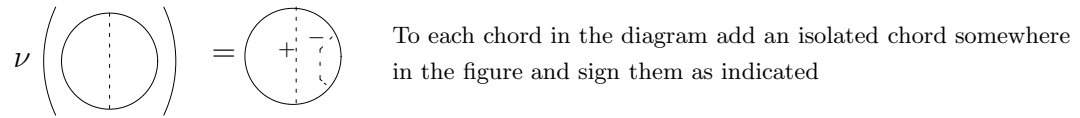
\[ 
  \printname{figure70c}
  \setlength{\unitlength}{0.5\standardunitlength}
  \begin{array}{c}  \hspace{-1.7mm}
    \raisebox{-8pt}{%
\begingroup\makeatletter\ifx\SetFigFont\undefined
\def\x#1#2#3#4#5#6#7\relax{\def\x{#1#2#3#4#5#6}}%
\expandafter\x\fmtname xxxxxx\relax \def\y{splain}%
\ifx\x\y   
\gdef\SetFigFont#1#2#3{%
  \ifnum #1<17\tiny\else \ifnum #1<20\small\else
  \ifnum #1<24\normalsize\else \ifnum #1<29\large\else
  \ifnum #1<34\Large\else \ifnum #1<41\LARGE\else
     \huge\fi\fi\fi\fi\fi\fi
  \csname #3\endcsname}%
\else
\gdef\SetFigFont#1#2#3{\begingroup
  \count@#1\relax \ifnum 25<\count@\count@25\fi
  \def\x{\endgroup\@setsize\SetFigFont{#2pt}}%
  \expandafter\x
    \csname \romannumeral\the\count@ pt\expandafter\endcsname
    \csname @\romannumeral\the\count@ pt\endcsname
  \csname #3\endcsname}%
\fi
\fi\endgroup
{\renewcommand{\dashlinestretch}{30}
\begin{picture}(11305,1385)(0,-10)
\put(431.250,683.000){\arc{3187.500}{5.8458}{6.7205}}
\put(1968.750,683.000){\arc{3187.500}{2.7043}{3.5789}}
\put(1200,683){\ellipse{1200}{1200}}
\put(3487,744){\ellipse{1200}{1200}}
\dashline{60.000}(1200,1283)(1200,83)
\dashline{60.000}(3525,1358)(3525,158)
\dashline{60.000}(3975,983)(3825,833)(3825,458)(3900,383)
\put(2475,683){\makebox(0,0)[lb]{\smash{{\mathmode{= }}}}}
\put(4725,908){\makebox(0,0)[lb]{\smash{{{\rm\scriptsize To each chord in the diagram add an isolated chord somewhere}}}}}
\put(0,608){\makebox(0,0)[lb]{\smash{{\mathmode{\nu}}}}}
\put(3300,758){\makebox(0,0)[lb]{\smash{{\mathmode{\scriptstyle +}}}}}
\put(3675,908){\makebox(0,0)[lb]{\smash{{\mathmode{\scriptstyle -}}}}}
\put(4725,458){\makebox(0,0)[lb]{\smash{{{\rm\scriptsize in the figure and sign them as indicated }}}}}
\put(10875,683){\makebox(0,0)[lb]{\smash{{\mathmode{\relax}}}}}
\end{picture}
}
 }
    \hspace{-1.9mm}
  \end{array}
 \]
Applying $\nu$ on the $4$-T relation yields the following two equalities
that assures that $\nu$ is a well defined map from $\mathcal{A}^{c}$
into $\mathcal{A}^{DD}$.
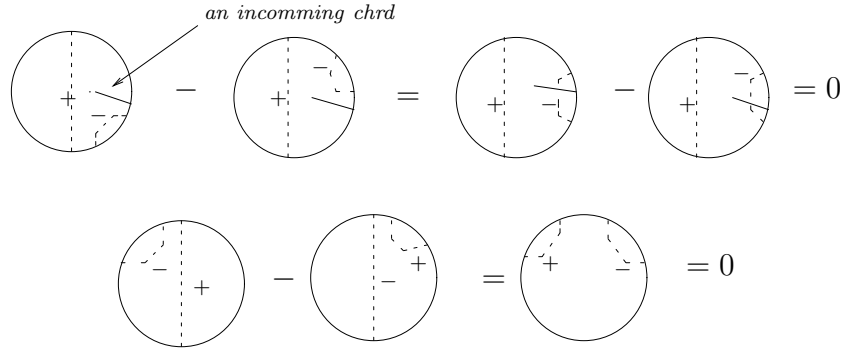
\[ 
  \printname{figure91}
  \setlength{\unitlength}{0.5\standardunitlength}
  \begin{array}{c}  \hspace{-1.7mm}
    \raisebox{-8pt}{\input draws/figure91.tex }
    \hspace{-1.9mm}
  \end{array}
 \]
\[ 
  \printname{figure92}
  \setlength{\unitlength}{0.5\standardunitlength}
  \begin{array}{c}  \hspace{-1.7mm}
    \raisebox{-8pt}{%
\begingroup\makeatletter\ifx\SetFigFont\undefined
\def\x#1#2#3#4#5#6#7\relax{\def\x{#1#2#3#4#5#6}}%
\expandafter\x\fmtname xxxxxx\relax \def\y{splain}%
\ifx\x\y   
\gdef\SetFigFont#1#2#3{%
  \ifnum #1<17\tiny\else \ifnum #1<20\small\else
  \ifnum #1<24\normalsize\else \ifnum #1<29\large\else
  \ifnum #1<34\Large\else \ifnum #1<41\LARGE\else
     \huge\fi\fi\fi\fi\fi\fi
  \csname #3\endcsname}%
\else
\gdef\SetFigFont#1#2#3{\begingroup
  \count@#1\relax \ifnum 25<\count@\count@25\fi
  \def\x{\endgroup\@setsize\SetFigFont{#2pt}}%
  \expandafter\x
    \csname \romannumeral\the\count@ pt\expandafter\endcsname
    \csname @\romannumeral\the\count@ pt\endcsname
  \csname #3\endcsname}%
\fi
\fi\endgroup
{\renewcommand{\dashlinestretch}{30}
\begin{picture}(7387,1689)(0,-10)
\put(1920,762){\makebox(0,0)[lb]{\smash{{\mathmode{-}}}}}
\put(6195,912){\makebox(0,0)[lb]{\smash{{\mathmode{\scriptstyle -}}}}}
\put(3270,762){\makebox(0,0)[lb]{\smash{{\mathmode{\scriptstyle -}}}}}
\put(4545,762){\makebox(0,0)[lb]{\smash{{\mathmode{=}}}}}
\put(7095,912){\makebox(0,0)[lb]{\smash{{\mathmode{= 0 }}}}}
\put(945,687){\makebox(0,0)[lb]{\smash{{\mathmode{\scriptstyle +}}}}}
\put(3645,987){\makebox(0,0)[lb]{\smash{{\mathmode{\scriptstyle +}}}}}
\put(5295,987){\makebox(0,0)[lb]{\smash{{\mathmode{\scriptstyle +}}}}}
\put(420,912){\makebox(0,0)[lb]{\smash{{\mathmode{\scriptstyle -}}}}}
\dashline{60.000}(795,1587)(795,12)
\dashline{60.000}(3195,1662)(3195,87)
\dashline{60.000}(5520,1587)(5520,1437)(5295,1137)(5070,1137)
\dashline{60.000}(6120,1587)(6120,1362)(6345,1062)(6570,1062)
\dashline{60.000}(570,1512)(570,1287)(345,1062)(45,1062)
\dashline{60.000}(3420,1587)(3420,1362)(3570,1212)(3870,1287)
\put(5820,875){\ellipse{1574}{1574}}
\put(3195,875){\ellipse{1574}{1574}}
\put(795,800){\ellipse{1574}{1574}}
\end{picture}
}
 }
    \hspace{-1.9mm}
  \end{array}
 \]
Notice that if D is a chord diagram then the following holds.
\begin{displaymath}
\iota \circ \nu (D) = D + \textrm{sum of diagrams each of which has an
                            isolated chord}
\end{displaymath}
As a corollary of the previous equality, we have that in the case
of unframed
knots the restriction of $\iota$ to $\mathcal{A}^{DD}$ is a 1-1
map into $\mathcal{A}^{c,r}$ which is the same as
$\mathcal{A}^{c}$ divided by the framing independence relation
which says that is a chord diagram has an isolated chord then it
is equal to zero. From the definition of $\nu$ it is clear that
$\iota$  is onto $\mathcal{A}^{c,r}$. The meaning of the corollary
is that in the case of unframed knots, the theory of finite type
invariant for unframed knots with a fixed self linking number is
the
same as the theory of Vassiliev invariants of framed knots.
\subsection{2-Component String Links}
Recall~\cite{ba:2} for the definition of a string link with m
components.
In the case of $2$-components string tangles we can construct a
map $\nu$ with the same behaviour as for knots. For every chord in
a chord diagram D, we add another chord, at the bottom of the
diagram and we put $+$ sign on the
original chord, and a $-$ sign on the added new chord. See the
following example.
\[ 
  \printname{iottan}
  \setlength{\unitlength}{0.5\standardunitlength}
  \begin{array}{c}  \hspace{-1.7mm}
    \raisebox{-8pt}{\input draws/iottan.tex }
    \hspace{-1.9mm}
  \end{array}
 \]
Using the relations as before shows that this map is well-defined.
As a corollary we have, for two components tangles as well, that
both theories, of Vassiliev and finite type invariants are the
same.
\subsubsection{Remark}
\begin{figure}
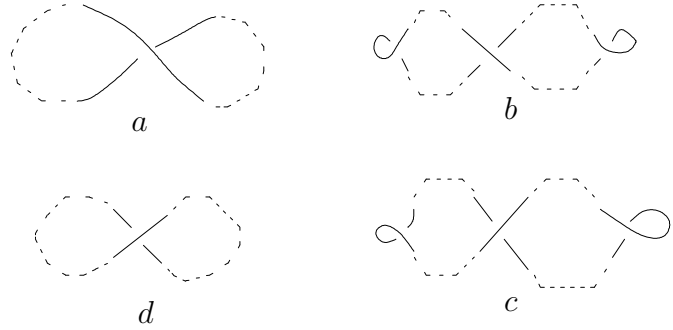

\[ 
  \printname{kink}
  \setlength{\unitlength}{0.5\standardunitlength}
  \begin{array}{c}  \hspace{-1.7mm}
    \raisebox{-8pt}{\input draws/kink.tex }
    \hspace{-1.9mm}
  \end{array}
 \]
\caption{A double flip(b to c) representing a single flip (a to d)}
\label{s:kink move}
\end{figure}

If we try to define a map similar to $\nu$ also for $m$-string tangles
where $m > 2$, this will not do, since the map will not be well-defined.
The order in which we add the negative chords is not irrelevant in this case.
We can interpret these corollaries geometrically as follows. If we
have a singular knot and a given intersection point on it, we can
artificially create two kinks, having signs opposite to each other
and as a result one of them has a sign opposite to the crossing we
wish to change when we resolve the given intersection point.
Performing a double-crossing change in which, the flipped
crossings are the crossing given by the intersection point and the
crossing which creates the kink with the opposite sign, is
topologically equivalent to performing the crossing flip of the
intersection point. The self linking number remains fixed during
this procedure, but still we should bear in mind the fact that the
framing is changed.
 This is indicated in figure ~\ref{s:kink move}.
For $2$-string tangles we have a similar way of presenting a single-crossing
change using a double-crossing change. Instead of kinks we create a trivial pair
of crossings for each crossing we wish to flip.\\
\subsubsection{The invariant $\mu_{ijk}$}
Calculating $\iota$ on the generator of the space of diagrams of
degree $2$ with $3$ components gives,
\[ 
  \printname{figure60}
  \setlength{\unitlength}{0.5\standardunitlength}
  \begin{array}{c}  \hspace{-1.7mm}
    \raisebox{-8pt}{%
\begingroup\makeatletter\ifx\SetFigFont\undefined
\def\x#1#2#3#4#5#6#7\relax{\def\x{#1#2#3#4#5#6}}%
\expandafter\x\fmtname xxxxxx\relax \def\y{splain}%
\ifx\x\y   
\gdef\SetFigFont#1#2#3{%
  \ifnum #1<17\tiny\else \ifnum #1<20\small\else
  \ifnum #1<24\normalsize\else \ifnum #1<29\large\else
  \ifnum #1<34\Large\else \ifnum #1<41\LARGE\else
     \huge\fi\fi\fi\fi\fi\fi
  \csname #3\endcsname}%
\else
\gdef\SetFigFont#1#2#3{\begingroup
  \count@#1\relax \ifnum 25<\count@\count@25\fi
  \def\x{\endgroup\@setsize\SetFigFont{#2pt}}%
  \expandafter\x
    \csname \romannumeral\the\count@ pt\expandafter\endcsname
    \csname @\romannumeral\the\count@ pt\endcsname
  \csname #3\endcsname}%
\fi
\fi\endgroup
{\renewcommand{\dashlinestretch}{30}
\begin{picture}(7320,2519)(0,-10)
\put(4237.500,1245.500){\arc{6978.628}{2.7790}{3.5041}}
\put(187.500,1245.500){\arc{6978.628}{5.9206}{6.6457}}
\put(6225,8){\makebox(0,0)[lb]{\smash{{\mathmode{k}}}}}
\put(6975,2333){\makebox(0,0)[lb]{\smash{{\mathmode{j}}}}}
\put(5475,2333){\makebox(0,0)[lb]{\smash{{\mathmode{i}}}}}
\put(3150,2408){\makebox(0,0)[lb]{\smash{{\mathmode{j}}}}}
\put(0,1358){\makebox(0,0)[lb]{\smash{{\mathmode{\iota}}}}}
\put(4275,1433){\makebox(0,0)[lb]{\smash{{\mathmode{=}}}}}
\put(1275,2333){\makebox(0,0)[lb]{\smash{{\mathmode{i}}}}}
\put(2100,8){\makebox(0,0)[lb]{\smash{{\mathmode{k}}}}}
\put(5587,1771){\ellipse{828}{828}}
\put(6900,1771){\ellipse{824}{824}}
\put(6225,721){\ellipse{824}{824}}
\dashline{60.000}(1800,1958)(2550,1958)
\dashline{60.000}(1725,1508)(2700,1508)
\dashline{60.000}(1875,1733)(2100,1583)(2250,1208)(2175,1133)
\put(2100,721){\ellipse{824}{824}}
\dashline{60.000}(6225,1433)(6525,1583)
\put(3000,1808){\ellipse{900}{900}}
\dashline{60.000}(6000,1583)(6225,1433)(6225,1133)
\dashline{60.000}(1575,1358)(1875,983)
\put(1425,1771){\ellipse{838}{838}}
\end{picture}
}
 }
    \hspace{-1.9mm}
  \end{array}
 \]
By~\cite{ha-li:1} the coefficient of the Chinese character
obtained from this diagram by removing the circles, in a certain
algebra structure, gives the weight system of the Milnor invariant
$\mu_{ijk}$.

\section{Open Questions}

We will now present some questions yet left open.
\begin{enumerate}
\item Are there any other relations in the spaces 
$\mathcal{D}^{DD}, \mathcal{D}^{\bigwedge}$ not presented here ?
\item What are the dimensions of $\mathcal{A}^{DD}, \mathcal{A}^{\bigwedge}$ in 
degrees higher than $2$ ? \item What is the kernel of $iota$ ?
\item How  should the map $\nu$ from $\mathcal{A}^c$ to
$\mathcal{A}^{DD}, \mathcal{A}^{\bigwedge}$ be defined in the general case, can 
it be well defined ?

\end{enumerate}

\ifx\undefined\bysame
        \newcommand{\bysame}{\leavevmode\hbox to3em{\hrulefill}\,}
\fi

\appendix
\section{Guided Tour Through Double Dating Diagrams}
In this section we use the relations in order to divide all DD
diagrams into their equivalence classes. Each of the equalities
listed here should be understood as equivalence modulo the
relations given in section $3$. Namely, two diagrams are
equivalent if the difference between them lies in the subspace
spanned by all the relations. For each equality we specify the
relation by which the equality is derived. When we say that some
diagram, D, is equal to zero we mean that either it has an
embedding $K_D$ on which all finite type invariant vanish or that
D is equivalent, via the relations to such a diagram.
 In the following we will refer to such $k_D$ as a geometric representative. All
equalities can be interpreted geometrically using appropriate
embeddings and drawing the crossing changes insinuated by the
relations. We will start with diagrams of degree $1$, and then
continue with degree-$2$ diagrams. Degree-$2$ diagrams are divided
into two groups. One group consists of diagrams with up to two
components.
 The second group consists of diagrams with three components.
 Equivalence classes of diagrams of degree $\leq$ 2 and more than $3$
components, are easily derived from the equalities given here. The
last element in the followin list should be understood as follows.
The diagram is
exactly one of the two summands that participate in the $2$-T
relation. The equality to one is for normalizing. We could have
chosen other non-zero elements as well, but since this diagram has
a close relation, which will be discussed in details at the end of
Appendix B, with the Milnor $\mu_{ijk}$ invariant which it self
takes the value one on three component- algebraically split links,
then we chose to use the same value here.
\subsubsection{Degree-$1$ Diagrams}
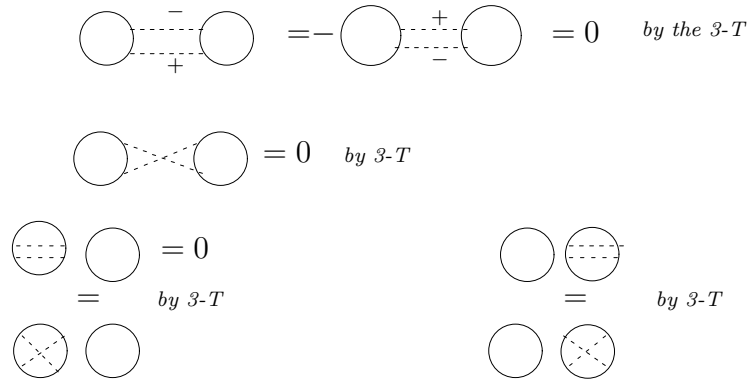
\[ 
  \printname{figure30}
  \setlength{\unitlength}{0.5\standardunitlength}
  \begin{array}{c}  \hspace{-1.7mm}
    \raisebox{-8pt}{\input draws/figure30.tex }
    \hspace{-1.9mm}
  \end{array}
 \]
\subsubsection{Degree-$2$ Diagrams}
Double dating diagrams of degree $2$ consist of two pairs of
chords attached to the circle components. There can be at most $3$
circles which are involved in a degree $2$ diagrams. A pair chords
is an isolated pair if
 there is an arc on each of the
circles this pair is attached to, so that no other chord is
attached to this arc. Otherwise the pair is non-isolated. The
complete list of degree $2$ double dating diagram shown below is
constructed by considering the various possibilities of having
isolated pairs. It is symmetric with regarding a permutation of
the set of the circle components of the diagram.
\begin{enumerate}
\item{ Two isolated pairs}
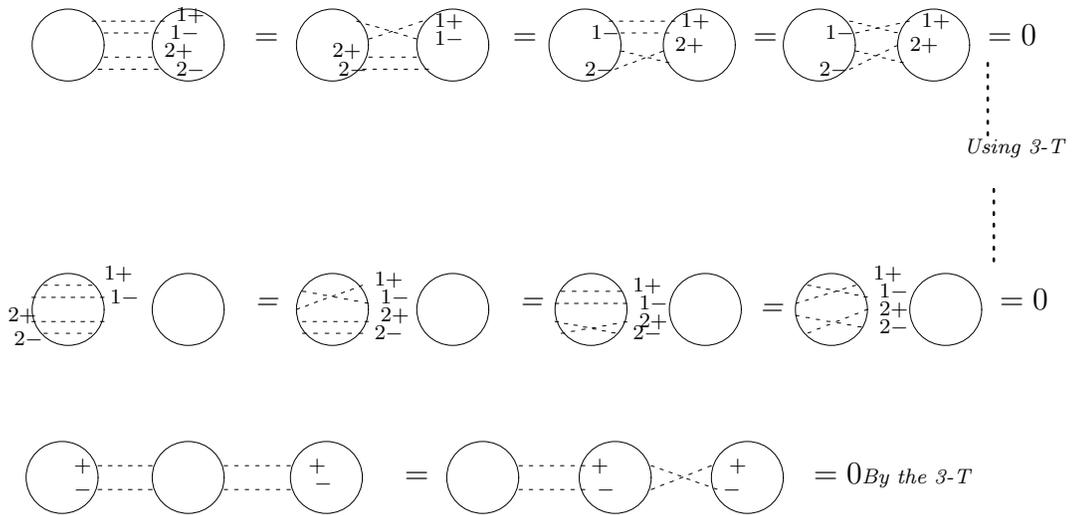
\[ 
  \printname{figure30aa}
  \setlength{\unitlength}{0.5\standardunitlength}
  \begin{array}{c}  \hspace{-1.7mm}
    \raisebox{-8pt}{\input draws/figure30aa.tex }
    \hspace{-1.9mm}
  \end{array}
 \]
\item{One isolated pair}
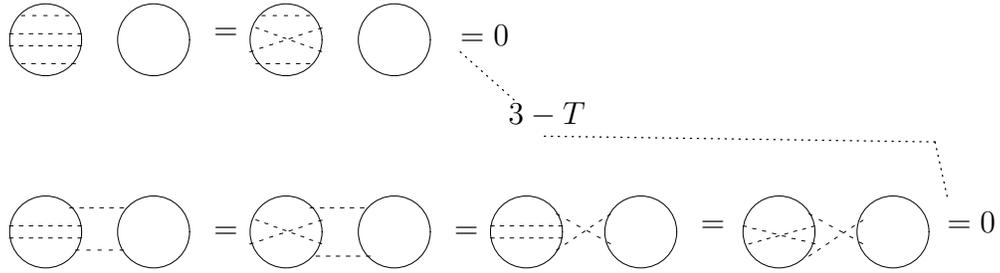
\[ 
  \printname{figure30ba}
  \setlength{\unitlength}{0.5\standardunitlength}
  \begin{array}{c}  \hspace{-1.7mm}
    \raisebox{-8pt}{\input draws/figure30ba.tex }
    \hspace{-1.9mm}
  \end{array}
 \]
\item{Two non-isolated pairs}
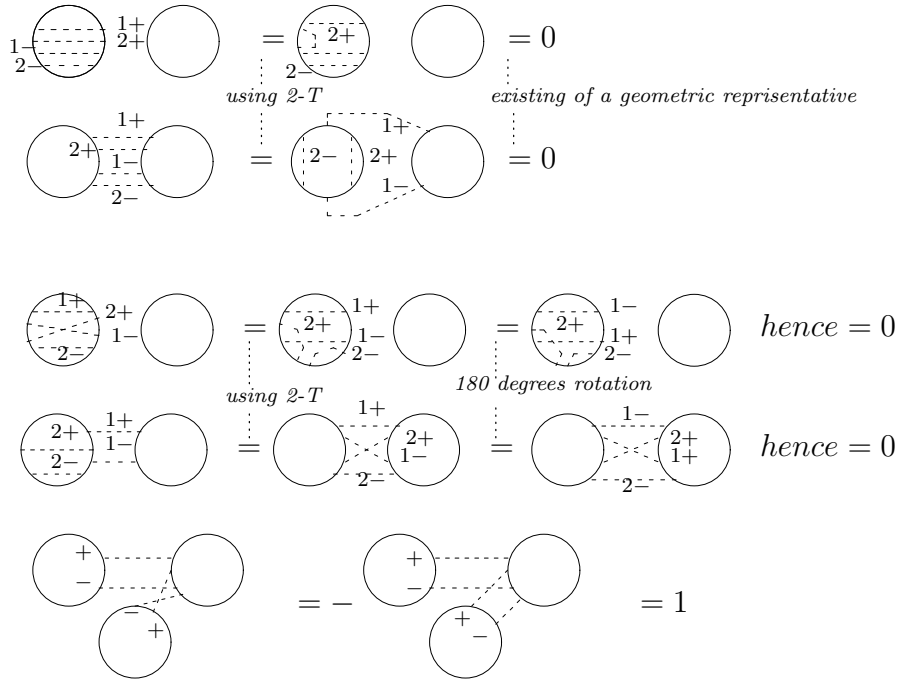
\[ 
  \printname{figure30ca}
  \setlength{\unitlength}{0.5\standardunitlength}
  \begin{array}{c}  \hspace{-1.7mm}
    \raisebox{-8pt}{\input draws/figure30ca.tex }
    \hspace{-1.9mm}
  \end{array}
 \]
\end{enumerate}
Now one can choose a representative for equivalence class and have
the result of theorem ~\ref{s:lowdim}
\section{Calculating Iota}
In this section we present the results obtained after applying $\iota$
to the relations we have for DD and $\bigwedge$ diagrams.
\[ 
  \printname{figurea1}
  \setlength{\unitlength}{0.5\standardunitlength}
  \begin{array}{c}  \hspace{-1.7mm}
    \raisebox{-8pt}{\input draws/figurea1.tex }
    \hspace{-1.9mm}
  \end{array}
 \]
\[ 
  \printname{figurea2}
  \setlength{\unitlength}{0.5\standardunitlength}
  \begin{array}{c}  \hspace{-1.7mm}
    \raisebox{-8pt}{\input draws/figurea2.tex }
    \hspace{-1.9mm}
  \end{array}
 \]
\[ 
  \printname{figurea3}
  \setlength{\unitlength}{0.5\standardunitlength}
  \begin{array}{c}  \hspace{-1.7mm}
    \raisebox{-8pt}{\input draws/figurea3.tex }
    \hspace{-1.9mm}
  \end{array}
 \]
Applying $\iota$ on the wood pecker relation gives the following
sum. Applying $\iota$ once more will yield a
sum of two 4-Ts and some terms that vanish in pairs.
\[ 
  \printname{figurea4}
  \setlength{\unitlength}{0.5\standardunitlength}
  \begin{array}{c}  \hspace{-1.7mm}
    \raisebox{-8pt}{\input draws/figurea4.tex }
    \hspace{-1.9mm}
  \end{array}
 \]

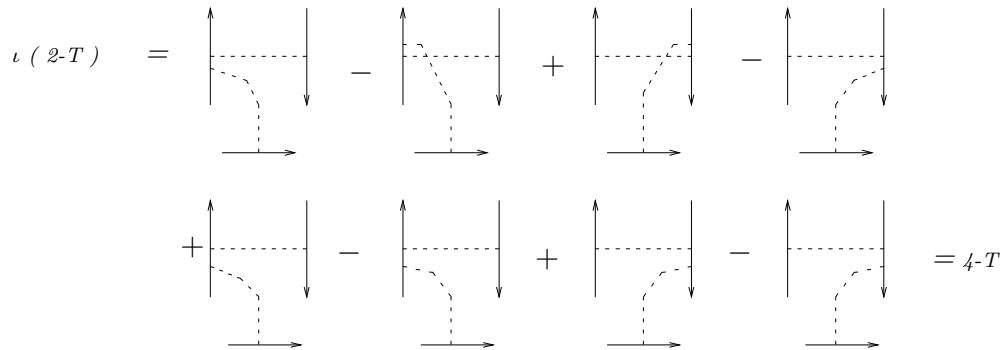
\end{document}